\documentclass[12pt]{amsart}
\usepackage{latexsym,amssymb}
\usepackage{amsfonts}
\usepackage{amsmath}
\usepackage{graphicx}
\usepackage[all]{xy}

\def\proof{\medskip\noindent{\sc Proof. }}
\def\EOP{\hfill$\Box$}%\RectangleBold}

\def\integ{{\mathbb Z}}

\def\be{\begin{equation}}
\def\ee{\end{equation}}

\newcommand{\fk}[1]{\mathfrak{#1}}
\newcommand{\mbf}[1]{\mathbf{#1}}

\newtheorem{theorem}{Theorem}[section]

\newtheorem{proposition}[theorem]{Proposition}

\newtheorem{remark}[theorem]{Remark}
\newtheorem{conjecture}[theorem]{Conjecture}
\newtheorem*{problem*}{Problem}

\theoremstyle{definition}
\newtheorem{definition/lemma}[theorem]{Definition/Lemma}

\newenvironment{changemargin}[2]{%
  \begin{list}{}{%
    \setlength{\topsep}{0pt}%
    \setlength{\leftmargin}{#1}%
    \setlength{\rightmargin}{#2}%
    \setlength{\listparindent}{\parindent}%
    \setlength{\itemindent}{\parindent}%
    \setlength{\parsep}{\parskip}%
  }%
  \item[]}{\end{list}}

\newcounter{ctr}

\begin{document}
\title[The toric ideal of a graphic matroid is generated by quadrics]
{The toric ideal of a graphic matroid \newline is generated by
quadrics}

\author{Jonah Blasiak}

\begin{abstract}
Describing minimal generating sets of toric ideals is a well-studied
and difficult problem. Neil White conjectured in 1980 that the toric
ideal associated to a matroid is generated by quadrics corresponding
to single element symmetric exchanges.  We give a combinatorial
proof of White's conjecture for graphic matroids.

\end{abstract}

\maketitle

\section{introduction}
Let $M$ be a matroid on the ground set $\{1, 2, \ldots n\}$. Fix a
field $k$ and define the polynomial ring $S_M$ to be $k[y_B: \;
\text{$B$ a base of $M$}]$.  Let $I_M$ be the kernel of the
$k$-algebra homomorphism $\theta_M: S_M \rightarrow k[x_1, \ldots,
x_n]$ that takes $y_B$ to $\prod_{i \in B} x_i$.  This is a toric
ideal as defined in \cite{Stu2}.

Given bases $B$ and $D$ of $M$, the well-known {\it symmetric
exchange property} states that for every $b \in B$ there exists a $d
\in D$ such that $B \cup d - b$ and $D \cup b - d$ are bases.  We
say that $b \in B$ {\it double swaps} into $D$; if $B \cup d - b$
and $D \cup b - d$ are bases, we say that $b \in B$ and $d \in D$
double swap.  Neil White made a conjecture in \cite{Wh2} about an
equivalence relation defined by certain symmetric exchange
properties and we state an algebraic reformulation.
\begin{conjecture}\label{whiteconjecture}
For any matroid $M$, the toric ideal $I_M$ is generated by the
quadratic binomials $y_{B_1} y_{B_2} - y_{D_1} y_{D_2}$ such that
the pair of bases $D_1, D_2$ can be obtained from the pair $B_1,
B_2$ by a double swap.
\end{conjecture}

The cycle matroid of a graph $G$, which we denote by $M(G)$, is the
matroid on the ground set $E(G)$ with a base for each spanning
forest of $G$.  A matroid is said to be {\it graphic} if it is the
cycle matroid of some graph.  We prove White's conjecture for
graphic matroids.
\begin{theorem}\label{graphictheorem}
If $M$ is a graphic matroid, then the toric ideal $I_M$ is generated
by the quadratic binomials $y_{B_1} y_{B_2} - y_{D_1} y_{D_2}$ such
that the pair of bases $D_1, D_2$ can be obtained from the pair
$B_1, B_2$ by a double swap.
\end{theorem}

To study $I_M$ in the context of toric ideals we need some notation.
Let $b$ be the number of bases of $M$ and let $A$ be the $n \times
b$ matrix whose columns are the zero-one incidence vectors of the
bases of $M$. The difference of two monomials is a {\it binomial}.
Given $\mathbf{u} \in \integ^b$ define $\mathbf{u_{+}}$
($\mathbf{u_{-}}$ resp.) to be $\mathbf{u}$ ($-\mathbf{u}$ resp.)
with negative coordinates replaced by zeros; we then have
$\mathbf{u} = \mathbf{u_+} - \mathbf{u_-}$.  The ideal $I_M$ is
spanned as a $k$-vector space by the binomials $\mathbf{y^{u_+}} -
\mathbf{y^{u_-}}$, where $\mathbf{u}$ runs over all integer vectors
in the kernel of $A$ \cite{Stu2}.  Ideals of this type, that is,
ideals generated by binomials $\mathbf{y^{u_+}} - \mathbf{y^{u_-}}$,
where $\mathbf{u}$ runs over integer vectors in the kernel of an
integer matrix, are toric ideals.  The set of vectors in $k^b$ that
vanish on all polynomials in a toric ideal is an affine toric
variety.  For each matroid $M$, the toric ideal $I_M$ is homogeneous
because every base has the same number of elements. Therefore $I_M$
(or any homogeneous toric ideal) defines a projective toric variety
$Y_M$ in $ k \mathbb{P}^{b-1}$  \cite{Stu2}.

White proves in \cite{Wh1} that
\begin{theorem}\label{projectivelynormal}
For any matroid $M$, the toric variety $Y_{M}$ is projectively
normal.
\end{theorem}
The following is a general conjecture about projectively normal
toric varieties \cite{Stu2}.
\begin{conjecture}\label{normal conjecture}
If the toric ideal $I$ defines a projectively normal $r$-dimensional
toric variety, then $I$ has a Gr\"obner basis consisting of
binomials of degree at most $r$.
\end{conjecture}
This conjecture restricted to toric varieties coming from matroids
neither implies nor is implied by White's conjecture. However it is
natural to ask whether the following variant of White's conjecture
holds (see, for instance, chapter 14 of \cite{Stu1} and \cite{HH})
and this does imply Conjecture \ref{normal conjecture} for toric
varieties $Y_M$ coming from a matroid $M$.
\begin{conjecture}\label{Groebner conjecture}
For any matroid $M$, the toric ideal $I_M$ has a Gr\"obner basis
consisting of quadratic binomials.
\end{conjecture}

%If a quadratic Gr\"obner basis is known, it allows for easy
%computation of many properties of the matroid polytope including
%volume [cite ???].  This volume determines the degree of the variety
%\cite{Stu2}.

White's conjecture can be posed as two separate conjectures.  The
following are both still open and together imply White's conjecture.
\begin{conjecture}\label{whiteconjecture A}
For any matroid $M$, the toric ideal $I_M$ is generated by quadratic
binomials.
\end{conjecture}
\begin{conjecture}\label{whiteconjecture B}
For any matroid $M$, the quadratic binomials of $I_M$ are in the
ideal generated by the binomials $y_{B_1} y_{B_2} - y_{D_1} y_{D_2}$
such that the pair of bases $D_1, D_2$ can be obtained from the pair
$B_1, B_2$ by a double swap.
\end{conjecture}
%watch out, does this formulation always allow permuting bases? yes

Sturmfels shows in chapter 14 of \cite{Stu1} that Conjecture
\ref{Groebner conjecture} holds for uniform matroids. One may also
ask the same questions about toric ideals coming from polymatroids.
Conca proves Conjecture \ref{whiteconjecture A}
for transversal polymatroids \cite{Con}.  %these matroids are products of rank 1 matroids I think
Caviglia, Elizalde, and Garc\'{i}a prove that both White's
conjecture and Conjecture \ref{Groebner conjecture} hold for a
certain class of polymatroids they call staircase polymatroids
\cite{CEG}.

In Section 2 we show that Conjecture \ref{whiteconjecture A} holds
for graphic matroids if certain graphs $\mathfrak{G}_k(M)$, defined
for $k \geq 3$, are connected for all graphic matroids $M$.
Similarly, Conjecture \ref{whiteconjecture B} holds for graphic
matroids if the graphs $\mathfrak{G}(M)$ are connected for all
graphic matroids $M$. In Section 3 we prove that the graphs
$\mathfrak{G}_k(M)$ are connected for any graphic matroid $M$.  In
Section 4 we prove that the graph $\mathfrak{G}(M)$ is connected for
any graphic matroid $M$.  In section 5 we discuss the difficulties
of extending our results to general matroids and pose some questions
along these lines.

\section{Reduction}
We show that the algebraic formulation of White's conjecture is
implied by a combinatorial condition similar to White's original
formulation.

Let $M$ be a matroid on a ground set of size $r(M) k$, where $r(M)$
denotes the rank of $M$. The {\it $k$-base graph} of $M$, which we
denote by $\mathfrak{G}_k(M)$, has as its vertex set the set of all
sets of $k$ disjoint bases (this is equivalent to the condition that
the union of the $k$ bases is the entire ground set). There is an
edge between $\{B_1,\dots,B_k\}$ and $\{D_1,\dots,D_k\}$ if and only
if $B_i = D_j$ for some $i,j$. We prove that Conjecture
\ref{whiteconjecture A} is implied by the connectivity of the
$k$-base graphs.  We prove the following proposition for a general
class of matroids $\mathfrak{C}$ that is closed under deletions and
adding parallel elements, but we will only apply this to the case
where $\mathfrak{C}$ is the set of graphic matroids.

\begin{proposition} \label{reduction}
Let $\mathfrak{C}$ be a collection of matroids that is closed under
deletions and adding parallel elements.  Suppose that for each $k
\geq 3$ and for every matroid $M$ in $\mathfrak{C}$ on a ground set
of size $r(M) k$ the $k$-base graph of $M$ is connected.  Then for
every matroid $M$ in $\mathfrak{C}$, $I_M$ is generated by quadratic
binomials.
\end{proposition}

\proof We will prove by induction on $k$ the statement that for
every $M \in \fk{C}$ and every binomial $b \in I_M$ of degree $k$,
$b$ is in the ideal generated by the quadrics of $I_M$.  This will
prove the proposition because, as mentioned in the introduction,
$I_M$ is spanned as a $k$-vector space by binomials.  For the base
case $k=2$ there is nothing to prove.  Suppose $k \geq 3$, $M$ is a
matroid in $\fk{C}$ on the ground set $\{1, 2,\ldots, n\}$, and $b$
is a binomial in $I_M$.  The binomial $b$ is necessarily of the form
$b = \prod_{i=1}^k y_{B_i} - \prod_{i=1}^k y_{D_i}$ for some bases
$B_1, \ldots, B_k, D_1, \ldots, D_k$ of $M$ such that the $B_i$ and
$D_i$ have the same multiset union.  We will show that $b$ is in the
ideal generated by the degree $k-1$ binomials of $I_M$.  By
induction, the degree $k-1$ binomials are in the ideal generated by
the quadrics of $I_M$ so this will complete the proof.

Put $\mathbf{x^S} = \theta_M(\prod_{i=1}^k y_{B_i})$ and let
$\mathbf{S}_i$ denote the $i^{th}$ component of $\mathbf{S}$. Define
$M'$ to be the matroid obtained from $M$ by replacing $i$ with
$\mathbf{S}_i$ parallel copies of $i$ for each $i$ in $\{1, \ldots,
n\}$; interpret ``replacing by zero parallel copies'' to mean
deleting, that is, delete those $i$ for which $\mathbf{S}_i = 0$.
There is a natural map $\alpha$ from the ground set of $M'$ to the
ground set of $M$ that takes each of the parallel copies of $i$ to
$i$.  A subset $X$ of the ground set of $M'$ is independent in $M'$
if and only if $\alpha(X)$ is independent in $M$.  This induces a
$k$-algebra homomorphism $\alpha_* : S_{M'} \rightarrow S_M$ defined
by $\alpha_*(y_B) = y_{\alpha(B)}$ for every base $B$ of $M'$.

Because the collection $\fk{C}$ is closed under deletions and adding
parallel elements, $M \in \fk{C}$ implies $M' \in \fk{C}$.  $M'$ has
a ground set of size $r(M') k$, and by assumption, the $k$-base
graph of $M'$ is connected.  Let $\mathbf{u_B}$ be a vertex of
$\fk{G}_k(M')$ such that $\alpha(\mathbf{u_B}) = \{B_1, \ldots,
B_k\}$ (here $\alpha$ is the natural extension of $\alpha$ to sets
of subsets of the ground of $M'$: $\alpha(\mbf{u_B}) = \{\alpha(X) |
X \in \mbf{u_B}\})$.  Such a $\mbf{u_B}$ exists by construction of
$M'$: simply split up the parallel copies of $i$, giving one to each
base in $\{B_1, \ldots, B_k\}$ containing $i$.  Let $\mathbf{u_D}$
be a vertex of $\fk{G}_k(M)$ such that $\alpha(\mathbf{u_D}) =
\{D_1, \ldots, D_k\}$.  Let $\mathbf{y^{u}} = \prod_{X \in \mbf{u}}
y_{X}$, as is customary when $\mbf{u}$ is identified with its
zero-one incidence vector.  Let $\mbf{u_0}, \mbf{u_1}, \ldots
,\mbf{u_t}$ be the vertices of a path between $\mbf{u_B} =
\mbf{u_0}$ and $\mbf{u_D} = \mbf{u_t}$ in $\fk{G}_k(M')$.  Then we
have
$$
  \sum_{i=1}^{t} \mbf{y^{u_{i-1}}} - \mbf{y^{u_{i}}}=
  \mbf{y^{u_0}} - \mbf{y^{u_t}}$$
and applying the map $\alpha_*$ we obtain \be \label{eq}
  \sum_{i=1}^{t} \mbf{y^{\alpha(u_{i-1})}} - \mbf{y^{\alpha(u_{i})}} =
 \mbf{y^{\alpha(u_0)}} - \mbf{y^{\alpha(u_t)}} =
 \prod_{i=1}^k y_{B_i} - \prod_{i=1}^k y_{D_i} = b .
 \ee
For $i = 1, \ldots, t$ there is a base $X \in \mbf{u_{i-1}} \cap
\mbf{u_{i}}$ which implies $\alpha(X) \in \mbf{\alpha(u_{i-1})} \cap
\mbf{\alpha(u_{i})}$.  This shows that $y_{\alpha(X)}$ may be
factored out of the binomial $\mbf{y^{\alpha(u_{i-1})}} -
\mbf{y^{\alpha(u_{i})}}$, and therefore (\ref{eq}) shows that $b$ is
in the ideal generated by the degree $k-1$ binomials of $I_M$. \EOP

The reduction for Conjecture \ref{whiteconjecture B} is similar.
Suppose $M$ is a rank $r$ matroid on a ground set of size $2r$.  The
{\it single exchange graph} of $M$, which we denote by
$\mathfrak{G}(M)$, is the graph with vertex set the set of ordered
2-tuples of bases of $M$, $(B_1, B_2)$, such that $B_1$ and $B_2$
are disjoint.  There is an edge between $(B_1, B_2)$ and $(D_1,
D_2)$ if and only if $|B_1 \cap D_1| = r-1$, or equivalently, $(D_1,
D_2)$ can be obtained from $(B_1, B_2)$ by a double swap. The above
proposition can be easily modified to show that: if for every $M$ in
$\fk{C}$ with a ground set of size $2 r(M)$ the single exchange
graph of $M$ is connected, then Conjecture \ref{whiteconjecture B}
holds for all matroids in $\fk{C}$.

\begin{remark}
{\rm Showing that $\mathfrak{G}$ is connected actually shows
slightly more than Conjecture \ref{whiteconjecture B}. In
$\mathfrak{G}$, $(B_1, B_2)$ is not adjacent to $(B_2,B_1)$ for
ranks larger than 1, however $y_{B_1} y_{B_2} - y_{B_2} y_{B_1} = 0$
is (trivially) in the ideal generated by quadrics corresponding to
single double swaps.  The stronger statement we prove here was also
conjectured by White in \cite{Wh2}. }

\end{remark}

\section{Proof of the graphic case}
We introduce some notation that is used in the main proof.  Let $G$
be a graph.  $V(G)$ and $E(G)$ denote the vertex and edge sets of
$G$.  If $v, v' \in G$, we abuse notation slightly and say that {\it
$v$ is connected to $v'$} or {\it $v$ and $v'$ are connected} to
mean that $v$ and $v'$ are in the same component. $d(v)$ denotes the
degree of $v$. We use $-$ to denote set minus and sometimes write a
one element set as the element itself rather than the element with
braces around it.

The following theorem together with Proposition \ref{reduction}
implies that Conjecture \ref{whiteconjecture A} holds for graphic
matroids.

\begin{theorem}\label{mainthm}
Let $G$ be a graph with $k r$ edges, where $r$ is the rank of
$M(G)$. If $k \geq 3$, then the $k$-base graph
$\mathfrak{G}_k(M(G))$ is connected.
\end{theorem}

\proof We prove the theorem by induction on $r$.  If $r =1$,
$\mathfrak{G}_k(M(G))$ is nonempty only if $G$ has no loops. If $G$
has no loops, $\mathfrak{G}_k(M(G))$ is a single vertex, which of
course is connected.  Now suppose $r>1$.  First observe that we can
assume $G$ is connected. If not, we may write $M(G)$ as the direct
sum of $M(G_1)$ and $M(G_2)$, where $G_1$ and $G_2$ are unions of
connected components of $G$.  By the inductive hypothesis, the
$k$-base graphs of $M(G_1)$ and $M(G_2)$ are connected.  The result
follows for $M(G)$ by Proposition 5 of \cite{Wh2}.

% on ground sets $E_1$ and
%$E_2$ respectively. Suppose $\{B_1,\dots,B_k\}$ and
%$\{D_1,\dots,D_k\}$ are vertices in $\mathfrak{G}_k(M(G))$.  Then
%$\{B_1 \cap E_1 ,\dots,B_k \cap E_1\}$ and $\{D_1 \cap E_1
%,\dots,D_k \cap E_1\}$ are vertices of $\mathfrak{G}_k(M_1)$.  By
%induction they are connected by a path
A key observation is that $M(G)$ has a cocircuit of size $\leq
2k-1$.  This is not true in general matroids and this is the most
essential way the graphic hypothesis is used. The graph $G$ has a
vertex $v$ of degree $\leq 2k - 1$ because $G$ has $r+1$ vertices
and $k r$ edges, making the average vertex degree $2 k
\frac{r}{r+1}$. The vertex $v$ is fixed throughout the proof.  Let
$C$ be the set of edges leaving $v$ and let $N(v)$ be the neighbors
of $v$.

We say a vertex $\{B_1,\dots,B_k\}$ of $\mathfrak{G}_k(M(G))$ is
{\it balanced} if $|B_i \cap C| \leq 2$ for each $i$.    We first
show that each vertex of $\mathfrak{G}_k(M(G))$ is connected to a
balanced vertex. We then show that any two balanced vertices that
have the same intersections with $C$ are connected. This is the
heart of the proof and where the inductive step is used. Finally, we
show that any two balanced vertices are connected.  These facts are
proved in this order as statements (1), (2), and (3), and these are
enough to show that $\mathfrak{G}_k(M(G))$ is connected.

\medskip{\it \noindent
\emph{\textbf{(1)}} Any vertex $\{B_1,\dots,B_k\}$ of
$\mathfrak{G}_k(M(G))$ is connected to a balanced vertex.}
\medskip

\begin{figure}[tbp]
\begin{changemargin}{-20pt}{0pt}
\begin{center}
\includegraphics [scale=.9]{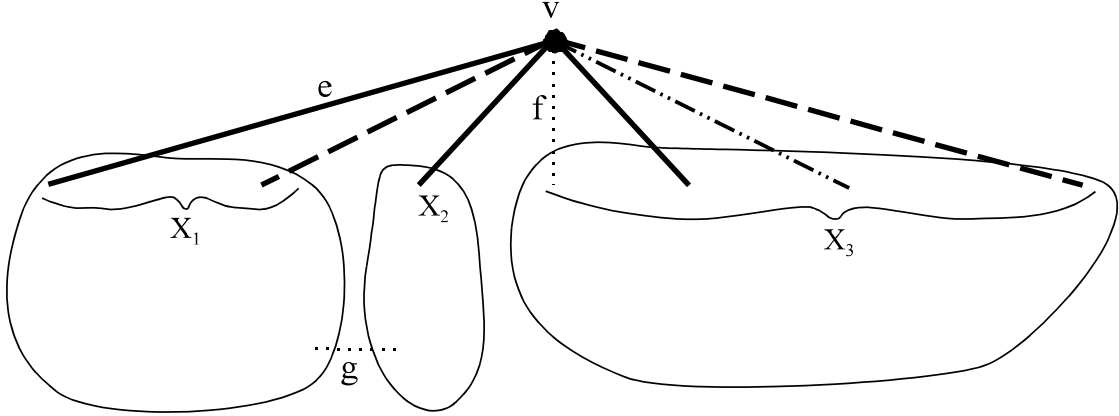}
\end{center}
\end{changemargin}
\caption{Each edge type corresponds to one of the bases $B_i$. The
normal edges correspond to $B_1$ and the dotted edges correspond to
$B_2$.  The blobs represent the components of $H$.}
\label{vertexstar}
\end{figure}

%begin proof
Let $\{B_1,\dots,B_k\}$ be a vertex of $\mathfrak{G}_k(M(G))$. Let
$S_i = B_i \cap C$.  Suppose that $\{B_1,\dots,B_k\}$ is not
balanced and (without loss of generality) $|S_1| > 2$.  Consider the
subgraph $H$ of $G$ with edge set $B_1 - C$ and vertex set $V(G) -
v$. It has $|S_1|$ components; the intersection of these components
with $N(v)$ partitions $N(v)$, and therefore $C$, into $|S_1|$
parts.  We denote this partition by $X_1 \cup \ldots \cup X_{|S_1|}
= C$.  See Figure \ref{vertexstar}. Note that $S_1$ intersects each
of the $X_i$ in size 1. As $d(v) \leq 2k-1$, without loss of
generality $|S_2| = 1$. Say $S_2 = \{f\}$ and $e \in S_1$ is an edge
not in the $X_i$ containing $f$ (This is $X_3$ in the figure, and $e
\in X_1$; all we need is that $e \notin X_3$). Now double swap $e$
out of $B_1$ and into $B_2$.  That is, there exists a $g \in B_2$
such that $B_1 \cup g - e$ and $B_2 \cup e - g$ are bases. The edge
$g$ is not in $C$ because $g \in B_2$, $B_2 \cap C = \{f\}$, and $f$
and $g$ are distinct; if $f = g$, then $B_1 \cup g - e$ intersects
$X_3$ in size 2 contradicting that it's a base. Therefore $|(B_1
\cup g - e) \cap C| = |S_1| - 1$ and $|(B_2 \cup e - g) \cap C| =
2$. By repeating such swaps we eventually obtain a balanced vertex.
This proves (1).

\medskip

Given a balanced vertex $\{B_1,\dots,B_k\}$, its {\it matching
graph} is the graph with vertex set $C$ and an edge with ends $B_i
\cap C$ for each $i$ such that $|B_i \cap C| = 2$.  Note that the
matching graph has vertices of degree at most one and at least one
isolated vertex.

\medskip{\it \noindent
\emph{\textbf{(2)}} If two balanced vertices $\{B_1,\dots,B_k\},
\{D_1,\dots,D_k\}$ of $\mathfrak{G}_k(M(G))$ have identical matching
graphs, then they are connected.}
\medskip

\begin{figure}
%\begin{changemargin}{-70pt}{0pt}
\includegraphics[scale=1.2,trim = 0 0 250 0, clip]{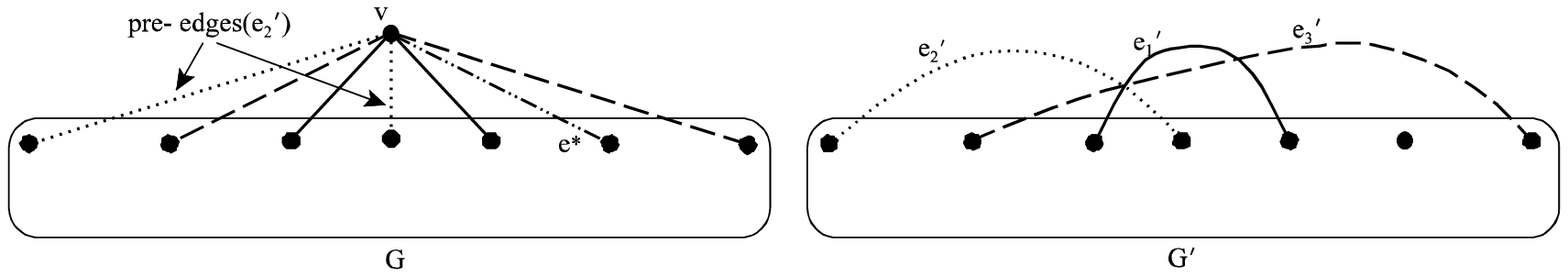}
\includegraphics[scale=1.2,trim = 250 0 0 0, clip]{Gprime.eps}
%\end{changemargin}
\caption{Each edge type corresponds to a base of $G$ and $G'$.}
\label{Gprime}
\end{figure}

We obtain a new graph $G'$ from $G$ as follows (see Figure
\ref{Gprime}): delete $v$ and for each $B_i$ with $|B_i \cap C| = 2$
add an edge between the vertices of $N(v)$ that are ends of the two
edges in $B_i \cap C$ (the subgraph of $G'$ induced by $N(v)$ is the
matching graph of $\{B_1,\dots,B_k\}$); call this new edge $e_i'$
and call $B_i \cap C$ the {\it pre-edges} of $e_i'$. Let $Z$ denote
the set of these special edges $e_i'$. Note that $|Z| \leq k-1$. For
$i$ such that $|B_i \cap C| = 2$, let $B_i' = B_i \cup e_i' - C$ and
$D_i' = D_i \cup e_i' - C$. Note that if we look at the subgraph of
$G$ with edge set $B_i$, then $B_i'$ is obtained by unsubdividing
$v$. For $i$ such that $|B_i \cap C| = 1$, put $B_i' = B_i - C$ and
$D_i' = D_i - C$. By induction on $r$, there is a path from
$\{B_1',\dots,B_k'\}$ to $\{D_1',\dots,D_k'\}$ in
$\mathfrak{G}_k(M(G'))$.  We will convert this to a path in
$\mathfrak{G}_k(M(G))$.

Given any set of $k$ disjoint bases of $M(G')$, we can reverse the
above process to produce $k$ disjoint bases of $M(G)$: if some base
$B$ of $M(G')$ intersects $Z$ in size $t > 0$, choose $t+1$ of the
pre-edges of $B \cap Z$ so that the resulting union with $B - Z$ is
a base of $M(G)$ (not all choices of $t+1$ edges will work, but at
least one will since $B-Z$ is a forest and $(B-Z)\cup
\text{pre-edges}(B\cap Z)$ spans $V(G)$). The $t-1$ pre-edges not
used will be added to bases not intersecting $Z$. Also, there exists
$e^* \in C$ that is not the pre-edge of any $e_i'$. This will always
be added to some base not intersecting $Z$.  We call this process of
taking a base of $M(G')$ and producing a base of $M(G)$ {\it pulling
back}, and we call the base of $M(G)$ the {\it pull back} of the
base of $M(G')$; we also use this terminology for sets of bases as
follows. To pull back a vertex $\{M_1',\dots,M_k'\}$ of
$\mathfrak{G}_k(M(G'))$, pull back the bases intersecting $Z$ first.
There are typically many choices for each of these pull backs, and
these choices can be made independently since the sets
$\text{pre-edges}(M_i' \cap Z)$ are disjoint.  Next, pull back the
bases not intersecting $Z$ by adding to each a single edge of $C$
not yet used by the other pull backs.

Pull each vertex in the path from $\{B_1',\dots,B_k'\}$ to
$\{D_1',\dots,D_k'\}$ back to a vertex of $\mathfrak{G}_k(M(G))$.
Now suppose that $\{M_1',\dots,M_k'\}$ and $\{N_1',\dots,N_k'\}$ are
consecutive vertices in the path in $\mathfrak{G}_k(M(G'))$. Without
loss of generality, $M_1' = N_1'$. Let $M_1$ and $N_1$ be the
corresponding pulled back bases of $M(G)$.  We want $M_1 = N_1$. If
$M_1'$ intersects $Z$, we can force $M_1 = N_1$ since the pull backs
of $M_1'$ and $N_1'$ do not depend on the pull backs of $M_i', N_i'$
for $i > 1$.

If $M_1'$ does not intersect $Z$, $M_1$ may differ from $N_1$ by one
element. Suppose $e^* \in M_j$. If $j \neq 1$, double swap $e^*$ of
$M_j$ with $M_1 \cap C$ of $M_1$ (this is possible because $|M_j
\cap C| = |M_1 \cap C| = 1$).  Denote the resulting set of bases by
$\{P_1,\dots,P_k\}$, and put $\{P_1,\dots,P_k\} = \{M_1,\dots,M_k\}$
in the case $j=1$. Do the same thing with $\{N_1,\dots,N_k\}$ and
$N_1$ to obtain $\{Q_1,\dots,Q_k\}$. The vertices
$\{P_1,\dots,P_k\}$ and $\{Q_1,\dots,Q_k\}$ are adjacent in
$\mathfrak{G}_k(M(G))$ because $P_1 = M_1' \cup e^* = N_1' \cup e^*
= Q_1$.  Therefore there is a path between $\{M_1,\dots,M_k\}$ and
$\{N_1,\dots,N_k\}$ in $\mathfrak{G}_k(M(G))$.  This proves that the
pulled back path can be patched up to make a path from
$\{B_1,\dots,B_k\}$ to $\{D_1,\dots,D_k\}$ in
$\mathfrak{G}_k(M(G))$.  This proves (2).

\bigskip

\smallskip{\it \noindent \emph{\textbf{(3)}} If $\{B_1,\dots,B_k\}$
and $\{D_1,\dots,D_k\}$ are balanced vertices of
$\mathfrak{G}_k(M(G))$, then there is a balanced vertex
$\{M_1,\dots,M_k\}$ connected to $\{B_1,\dots,B_k\}$ and a balanced
vertex $\{N_1,\dots,N_k\}$ connected to $\{D_1,\dots,D_k\}$ such
that $\{M_1,\dots,M_k\}$ and $\{N_1,\dots,N_k\}$ have the same
matching graph.}
\medskip

First note that (2) and (3) together show that any two balanced
vertices are connected, and therefore proving (3) will complete the
proof of Theorem \ref{mainthm}.  We prove (3) by rearranging the
parts of the bases that intersect $C$ without changing the other
parts.  Although the proof is rather involved, it is not hard to
convince oneself that the result is true by trying small values of
$k$.  Proving the result for $k=3$ and $d(v) = 4$ only requires
checking a few cases.

A {\it valid move} on a matching graph $H$ of a vertex
$\{B_1,\dots,B_k\}$ is a change in the matching graph from $H$ to
$H'$ such that there is a vertex connected to $\{B_1,\dots,B_k\}$
with matching graph $H'$. First we show the existence of certain
valid moves and then we show that these are enough to prove (3).

\smallskip{\it \noindent \emph{\textbf{(A)}} Suppose $\{B_1,\dots,B_k\}$ is a
balanced vertex with matching graph $H$, $(e_1, e_2)$ and $(e_3,
e_4)$ are edges of $H$, and $e_5$ is an isolated vertex.  Then at
least one of (i) and (ii) and (isomorphically) at least one of (iii)
and (iv) are valid moves on $H$.  Furthermore, if (v) and (vi) are
not valid moves, then either (i) and (ii) are both valid or (iii)
and (iv) are both valid. \setcounter{ctr}{0}
\begin{list}{\emph{(\roman{ctr})}} {\usecounter{ctr} \setlength{\itemsep}{1pt} \setlength{\topsep}{2pt}}
\item Deleting $(e_1, e_2)$ and adding $(e_1,e_5)$.
\item Deleting $(e_1, e_2)$ and adding $(e_2,e_5)$.
\item Deleting $(e_3, e_4)$ and adding $(e_3,e_5)$.
\item Deleting $(e_3, e_4)$ and adding $(e_4,e_5)$.
\item Deleting $\{(e_1, e_2), (e_3, e_4)\}$ and adding $\{(e_1, e_3), (e_2,
e_4)\}$.
\item Deleting $\{(e_1, e_2), (e_3, e_4)\}$ and adding $\{(e_1, e_4), (e_2,
e_3)\}$.
\end{list}
}
\smallskip

We work again with double swaps.  Suppose $B_1 \cap C = e_1 \cup
e_2$, $B_2 \cap C = e_3 \cup e_4$, and $B_3 \cap C = e_5$.  Recall
from the proof of (1) that $B_1$ and $B_2$ determine partitions of
$C$ into two parts.  Suppose that $B_1$ determines the partition $L
\cup R = C = V(H)$ and $B_2$ determines $T \cup B = C = V(H)$.  For
$X \subset C$, $X \cup (B_1 - C)$ is a base if and only if $X$
intersects $L$ and $R$ in size 1 (a similar statement holds for
$B_2$).  Now to show the first part of (A), double swap $e_5 \in
B_3$ into $B_1$. The edge $e_5$ must be swapped with something in
$C$, so either (i) or (ii) holds. The same argument shows (iii) or
(iv) holds.

\begin{figure}
\begin{center}
\includegraphics [scale=1]{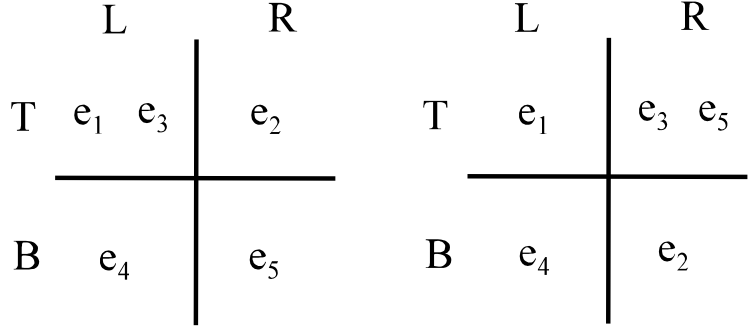}
\end{center}
\caption{Two possibilities for positions of $e_1, \ldots, e_5$ in
the partitions $L\cup R, T \cup B = V(H)$.  In the example on the
left, moves (i), (iii), and (vi) are valid.  In the example on the
right, moves (i),(ii), and (iv) are valid.} \label{partition}
\end{figure}

Consider the representation of the partitions $L \cup R = T \cup B =
V(H)$ as shown in Figure \ref{partition}.  The four regions
correspond to the sets $L \cap T, L \cap B, R \cap T$, and $R \cap
B$. If one of the regions contains two of $e_1, e_2, e_3, e_4$ (as
in the left example), then these elements can be double swapped.
This means we can replace $(e_1, e_2), (e_3, e_4)$ by either $(e_1,
e_3), (e_2, e_4)$ or $(e_1, e_4), (e_2, e_3)$ in the matching graph
and the resulting matching graph is realized by some vertex
connected to $\{B_1,\dots,B_k\}$; either (v) or (vi) is a valid
move.  If none of the regions contains two of $e_1, \ldots, e_4$,
then we are in a situation isomorphic to the right example of Figure
\ref{partition}. In this case (i) holds because $B_1 \cup e_5 - e_2$
and $B_3 \cup e_2 - e_5$ are bases, but in addition (ii) holds. This
is because $B_1 \cup \{e_3, e_4\}- \{e_1, e_2\}$, $B_2 \cup \{e_2,
e_5\}- \{e_3, e_4\}$, and $B_3 \cup e_1 - e_5$ are bases (if $e_5$
is in another region, (iii) and (iv) may hold instead of (i) and
(ii)). This proves (A).

\smallskip
\smallskip{\it \noindent \emph{\textbf{(B)}} Let $H$ and $H'$ be graphs
on the same vertex set both with maximum vertex degree 1 and the
same number, $t$, of isolated vertices, where $t > 0$.  It is
possible to get from $H$ to $H'$ by a sequence of valid moves of the
kind described in (A).}
\smallskip

We prove this by induction on $|V(H)|$.  The base case is when $H$
and $H'$ are both a single vertex.  Let $I_H$ be the set of vertices
that can be made isolated in $H$ after at most one valid move (two
vertices in this set don't have to be able to be isolated at the
same time). Define $I_{H'}$ similarly.  Using the moves (i) and
(ii), we see that $|I_H|, |I_{H'}| \geq t + |E(H)|$. Since $t +
|E(H)| > |V(H)|/2$, there is a vertex $x$ in $I_H \cap I_{H'}$. By
possibly redefining $H$ (or $H'$) to be a graph one move away from
it, we may assume that $x$ is isolated in $H$ and $H'$. If $t>1$,
delete $x$ from $H$ and $H'$, and the result follows by induction.

The case $t=1$ remains.  Consider the valid moves that make $x$ the
end of an edge: let $N_H$ be the set of vertices that can pair up
with $x$ after one move on $H$, and define $N_{H'}$ similarly.  We
have $|N_H| , |N_{H'}| \geq |E(H)|.$  If a move of type (i) and of
type (ii) are valid on $H$ then $|N_{H}| > |E(H)|$, and therefore
there exists $y \in (N_H \cap N_{H'})$.  Next, make the moves so
that both graphs have the common edge $(x,y)$.  Delete this edge
from both graphs and the result follows by induction.  For the rest
of the proof we may assume $x$ stays isolated and that for each edge
in $H$ and each edge in $H'$ the moves (i) and (ii) are not both
valid (we may also assume this for any graph we reach from $H$ or
$H'$ by a sequence of valid moves that keeps $x$ isolated). This
implies that for every pair of edges in $H$ and $H'$, either (v) or
(vi) is a valid move.

\begin{figure}[tbp]
\begin{center}
\begin{changemargin}{-50pt}{0pt}
\includegraphics [scale=1]{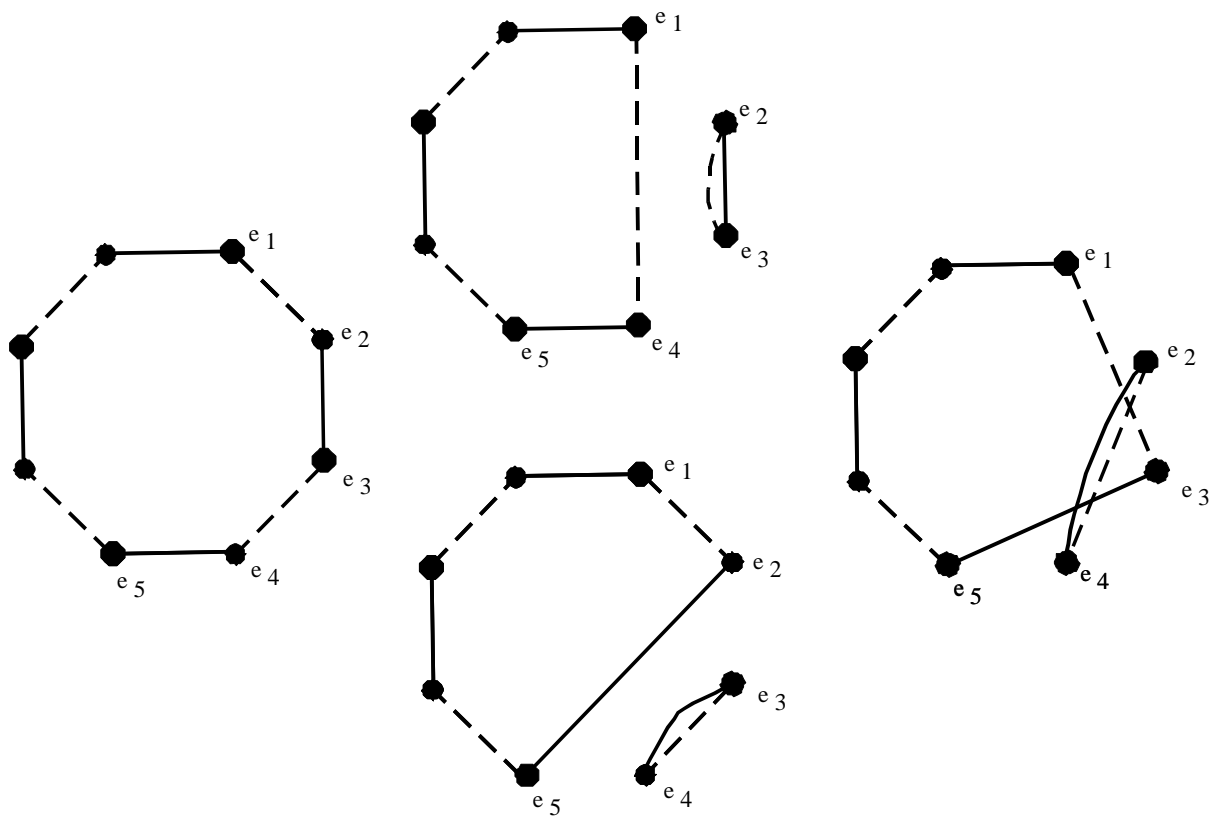}
\end{changemargin}
\end{center}
\caption{On the left is $J$, where dashed edges are edges of $H'$
and normal edges are edges of $H$.  The other graphs show $J$ after
possible moves on $H$ and $H'$ as described in the proof of (B).}
\label{8cycle}
\end{figure}

For the rest of proof, we modify the statement we are proving by
induction: we no longer require the graphs to have an isolated
vertex, but for each pair of edges either (v) or (vi) is valid. We
will prove this statement for the graphs $H-x$ and $H'-x$. Consider
the graph $J$ with vertex set $V(H)-x$ and edge set $E(H) \cup
E(H')$. It is 2-regular, and therefore a disjoint union of cycles.
If there is more than one component, split up $V(H)-x$ according to
the components and win by induction. The remaining case is if $J$ is
a cycle.  If $J$ is a 2-cycle, $H = H'$. If $|V(J)| \geq 4$, it
suffices to consider 5 consecutive vertices $e_1, e_2, e_3, e_4,
e_5$ as in Figure \ref{8cycle} ($e_1 = e_5$ if $J$ is a 4-cycle, but
the proof still works).  If replacing $(e_1, e_2)$ and $(e_3, e_4)$
by $(e_1, e_4)$ and $(e_2, e_3)$ is a valid move on $H'$, we get a
2-cycle and win by induction (as in the top graph of Figure
\ref{8cycle}). The same thing happens if replacing $(e_2, e_3)$ and
$(e_4, e_5)$ by $(e_2, e_5)$ and $(e_3, e_4)$ is a valid move on $H$
(as in the bottom graph). If neither of these is a valid move, then
(replacing $(e_1, e_2)$ and $(e_3, e_4)$ by $(e_1, e_3)$ and $(e_2,
e_4)$) and (replacing $(e_2, e_3)$ and $(e_4, e_5)$ by $(e_2, e_4)$
and $(e_3, e_5)$) are valid moves (as in the right graph). Delete
the ends of the common edge $(e_2, e_4)$ and win by induction.

\medskip

(B) implies (3) by letting $H$ be the matching graph of $\{B_1,
\ldots, B_k\}$ and $H'$ be the matching graph of $\{D_1, \ldots,
D_k\}$.  A sequence of valid moves beginning with $H$ yields a path
from $\{B_1,\dots,B_k\}$ to $\{M_1,\dots,M_k\}$ and a sequence of
valid moves beginning with $H'$ yields a path from
$\{D_1,\dots,D_k\}$ to $\{N_1,\dots,N_k\}$. (B) says that we can
find moves so that $\{M_1,\dots,M_k\}$ and $\{N_1,\dots,N_k\}$ have
the same matching graph.  \EOP

\section{Quadrics are generated by one element exchanges}
The following theorem together with the modified version of
Proposition \ref{reduction} shows that Conjecture
\ref{whiteconjecture B} holds for graphic matroids.  This will
complete the proof of Theorem \ref{graphictheorem}.

\begin{theorem}\label{frakGconnected}
Let $G$ be a graph with $2 r$ edges, where $r$ is the rank of
$M(G)$.  Then the single exchange graph $\mathfrak{G}(M(G))$ is
connected.
\end{theorem}

\proof The proof is very similar to the proof of Theorem
\ref{mainthm}.  We do induction on $r$.  We can assume $G$ is
connected for the same reason as before.  And again, we have that
there is a vertex $v$ of degree at most 3, which we fix throughout
the proof.  Let $C$ be the set of edges leaving $v$. There is no
need to balance the vertices because there is only one possibility
for the sizes of the intersections of two bases with $C$ (if $d(v) =
3$, one base intersects $C$ in size 2 and the other in size 1; if
$d(v)=2$ both bases intersect $C$ in size 1). As before, define the
matching graph of a vertex $(B_1, B_2)$ to be the graph with vertex
set $C$ and an edge with ends $B_i \cap C$ for $i$ such that $|B_i
\cap C| = 2$.  Note that the matching graph ignores the order of
$B_1$ and $B_2$; we are careful to remember this when proving (1)
below.

 We need to show that any two vertices of
$\mathfrak{G}$ that have the same matching graph are connected. This
is enough to prove the theorem since statement (3) of the proof of
Theorem \ref{mainthm} holds for $k=2$ with the same proof (although
a much simpler argument would do in this case).

\medskip{\it \noindent
\emph{\textbf{(1)}} If two vertices $(B_1,B_2), (D_1,D_2)$ of
$\mathfrak{G}(M(G))$ have the same matching graph, then they are
connected.}
\medskip

We obtain a new graph $G'$ from $G$ using the same construction as
in the proof of Theorem \ref{mainthm}. Let $B_1, B_2', D_1', D_2'$
be the bases of $M(G')$ defined in the proof of Theorem
\ref{mainthm}. By induction on $r$, there is a path from $(B_1',
B_2')$ to $(D_1', D_2')$ in $\mathfrak{G}(M(G'))$. We will convert
this to a path in $\mathfrak{G}(M(G))$.  Note that in this case the
pull backs are unique.

Pull each vertex in the path from $(B_1', B_2')$ to $(D_1',D_2')$
back to a vertex of $\mathfrak{G}(M(G))$. Now suppose that
$(M_1',M_2')$ and $(N_1',N_2')$ are consecutive vertices in the path
in $\mathfrak{G}(M(G'))$.  Let $M_1, M_2, N_1, N_2$ be the
corresponding pulled back bases of $M(G)$.  If $d(v)=2$, $(M_1,
M_2)$ is adjacent to $(N_1,N_2)$ and we are done.  If $d(v) = 3$,
observe that $Z$ consists of a single edge $e_1'$ and $C =
\text{pre-edges}(e_1') \cup e^*$. Without loss of generality, $e_1'
\in M_1'$. If $e_1' \in N_1'$, then $M_1$ and $N_1$ differ by only
one element and are therefore adjacent in $\mathfrak{G}(M(G))$.

If $e_1' \in N_2'$, then let $\{a,b\} =\text{pre-edges}(e_1')$.
Double swap $e^* \in M_2$ with $M_1$ to obtain a vertex $(P_1, P_2)$
adjacent to $(M_1, M_2)$. The edge $e^*$ must double swap with
something in $C$ (say, $a$), because otherwise $P_2$ would not
intersect $C$ contradicting that it's a base.  We know that $(N_2' -
e_1') \subset M_2'$ and we can rewrite this as $(N_2 - \{a,b\})
\subset (M_2 - e^*)$.  Add $\{a,b\}$ to both these sets to obtain
$N_2 \subset (P_2 \cup b)$ and therefore $|P_2 \cap N_2| = r-1$.
This shows that $(P_1, P_2)$ and $(N_1, N_2)$ are adjacent and thus
the pulled back path can be patched up to make a path from
$(B_1,B_2)$ to $(D_1, D_2)$ in $\mathfrak{G}(M(G))$. \EOP

\section{Future Work}
The proofs of Theorems \ref{mainthm} and \ref{frakGconnected} depend
heavily on the graphic assumption.  However, it seems possible to
convert many of the techniques to the general case.  For instance,
instead of choosing $C$ to be the edges leaving a vertex, we could
take $C$ to be a cocircuit.  There is an analog of the construction
of $G'$ for any cocircuit of a matroid. One thing that can
definitely not be generalized is the existence of a small cocircuit
and this is crucial to the proofs. For instance, there are uniform
matroids with arbitrarily large minimum cocircuit size for fixed
$k$.

Part $(3)$ of Theorem \ref{mainthm} at first seemed like a
digression from the main content of the proof and theorem, and a
fun, but not very significant, result on its own. However, the
analogous statement of $(3)$ for general cocircuits may actually be
rather deep. We will not state the exact generalization of $(3)$,
but it suggests the following question. Given matroids $M$ and $N$
on the ground set $E$ and $X \subset E$, define $r_{M \cap N}(X)$ to
be the maximum size of an independent set in $X$ common to $M$ and
$N$. Given matroids $M_1, \ldots, M_k$ and $N_1, \ldots, N_k$ all on
the ground set $E$, we define their {\it matching intersection rank}
to be $$ \max_{\pi \in S_k} \Bigl( \max_{X_1 \sqcup \ldots \sqcup
X_k = E} \sum_{i=1}^k r_{M_i \cap N_{\pi(i)}}(X_i) \Bigr).$$
%The
%quantity in parentheses is essentially the rank of $E$ in the union
%$(M_1 \cap N_{\pi(1)}) \cup \ldots \cup (M_k \cap N_{\pi(k)})$,
%although this union operation is typically only defined for
%matroids, and the rank function $r_{M \cap N}$ does not give the
%rank function of a matroid in general.
\medskip \noindent \textbf{Problem.}
 Suppose $\{B_1,\dots,B_k\}$ and $\{D_1,\dots,D_k\}$ are balanced
vertices of $\mathfrak{G}_k(M)$ with respect to some cocircuit $C$.
Here we take balanced to mean that the intersection sizes of the
bases with $C$ are less than one away from average intersection
size. Let $M_i = M.(B_i-C)|C$ and $N_i = M.(D_i-C)|C$, where $.$
denotes contraction and $|C$ means deleting everything not in $C$.
Determine conditions on $C$ under which the matching intersection
rank of $M_1, \ldots, M_k$ and $N_1, \ldots, N_k$ is $|C|$.
\medskip

This problem and the notion of matching intersection rank lead to
two general questions, but we have not been able to formulate
specific conjectures along these lines. Is there a generalization of
the matroid union and intersection theorems that says something
about matching intersection rank?  Does White's conjecture
generalize to a statement that involves bases of more than one
matroid?

\section{Acknowledgments}
The author thanks Bernd Sturmfels for many helpful conversations.


\begin{thebibliography}{9}
\bibitem{CEG} G. Caviglia and L. Garc\'{i}a, private communication, October 2005.

\bibitem{Con} A. Conca, {\it Linear spaces, transversal polymatroids and ASL
domains}, preprint 2005.

\bibitem{HH} J. Herzog and T. Hibi, {\it Discrete polymatroids
}, J. Algebraic Combin. {\bf 16} (2002), 239-268.

\bibitem{Stu1} B. Sturmfels, {\it Gr\"obner Bases and Convex
Polytopes}, American Mathematical Sociey, University Lecture Series,
Vol. 8, Providence, RI, 1995.

\bibitem{Stu2} B. Sturmfels, {\it Equations defining toric varieties},
Proc. Sympos. Pure Math. {\bf 62} (1997), 437-449.

\bibitem{Wh1} N. White, {\it The basis monomial ring of a matroid
}, Advances in Math. {\bf 24} (1977), 292-297.

\bibitem{Wh2} N. White, {\it A unique exchange property for bases},
 Linear Algebra and its Applications {\bf 31} (1980), 81-91.

\end{thebibliography}
\end{document}